\documentclass[11pt]{article}
\usepackage[T1]{fontenc}
\usepackage{lmodern}
\usepackage[letterpaper,textwidth=6.05in,textheight=9.25in,centering]{geometry}
\usepackage{amsmath,amssymb,amsthm}
\usepackage{microtype}
\usepackage{booktabs}
\usepackage{tikz}
\usepackage[hidelinks]{hyperref}

\newtheorem{theorem}{Theorem}[section]
\newtheorem{lemma}[theorem]{Lemma}
\newtheorem{proposition}[theorem]{Proposition}
\newcommand{\PSL}{\operatorname{PSL}_2(\mathbb C)}
\newcommand{\SL}{\operatorname{SL}_2(\mathbb C)}
\newcommand{\CP}{\mathbb{CP}^1}
\newcommand{\Vol}{\operatorname{Vol}}
\newcommand{\fix}{\operatorname{fp}}

\title{Every hyperbolic knot has an arc-faithful closed braid diagram}
\author{Haibin Hang\\[2pt]
  {\small\href{mailto:haibin.hang312@gmail.com}{haibin.hang312@gmail.com}}}
\date{}
\hypersetup{
  pdftitle={Every hyperbolic knot has an arc-faithful closed braid diagram},
  pdfauthor={Haibin Hang},
  pdfsubject={Braid words and non-pinched octahedral coordinates},
  pdfkeywords={hyperbolic knots, closed braids, octahedral decomposition}
}

\makeatletter
\renewcommand{\@maketitle}{%
  \begin{center}
    {\Large\bfseries \@title\par}
    \ifx\@author\@empty\else \vskip .7em {\large \@author\par}\fi
    \ifx\@date\@empty\else \vskip .4em \@date\par\fi
  \end{center}%
}
\makeatother

\newcommand{\braidlayer}[3]{%
  \pgfmathtruncatemacro{\rightpos}{#1+1}%
  \pgfmathtruncatemacro{\idlepos}{2-2*#1}%
  \pgfmathtruncatemacro{\lowerlevel}{#2-1}%
  \draw (\idlepos,#2) -- (\idlepos,\lowerlevel);
  \ifnum#3>0\relax
    \def\overstart{\rightpos}\def\overend{#1}%
  \else
    \def\overstart{#1}\def\overend{\rightpos}%
  \fi
  \draw (\overend,#2) .. controls +(0,-.43) and +(0,.43)
    .. (\overstart,\lowerlevel);
  \begin{scope}
    \clip ({#1+.5},{#2-.5}) circle[radius=3.3pt];
    \draw[white,line width=3.6pt] (\overstart,#2)
      .. controls +(0,-.43) and +(0,.43) .. (\overend,\lowerlevel);
  \end{scope}
  \draw (\overstart,#2) .. controls +(0,-.43) and +(0,.43)
    .. (\overend,\lowerlevel);
}

\begin{document}
\maketitle

\begin{abstract}
Every hyperbolic knot admits an arc-faithful closed braid diagram, answering
Conjecture~29 of McPhail--Snyder. More precisely, any braid representing the knot
can be reworded in the same braid group so that the complete holonomy has
distinct incoming parabolic fixed points at every crossing. The construction
replaces a pinched crossing by five non-pinched crossings when a suitable
third strand is adjacent. On \(n\) strands it uses at most \(2n-1\) letters
per replaced crossing. The resulting diagram realizes the complete holonomy
by nondegenerate shapes in its four-term octahedral decomposition.
\end{abstract}

\noindent{\small\textit{2020 Mathematics Subject Classification.}
57K32, 20F36, 57K10.}

\section{Introduction}\label{sec:intro}

Let \(M=S^3\setminus K\) be the complement of a hyperbolic knot \(K\).
A diagram \(D\) gives two descriptions of \(M\). The Wirtinger presentation
encodes \(\pi_1(M)\) by meridians and crossing relations. The octahedral decomposition
places one octahedron at each crossing and decomposes \(M\) with two
auxiliary points removed. In the four-term version each octahedron is
subdivided into four tetrahedra. The shape of a tetrahedron is a
cross-ratio of its developed ideal vertices. The gluing equations
assemble these shapes into diagram coordinates for the holonomy
\cite{KKY,McPhailSnyder}.

The complete structure of \(M\) has a discrete faithful holonomy
\(\rho_{\mathrm{hyp}}\colon\pi_1(M)\to\PSL\) with nontrivial parabolic
meridians. A crossing is \emph{pinched} if the incoming meridians have the
same holonomy fixed point. This coincidence forces shape parameters to take
the degenerate value \(1\). McPhail--Snyder calls \(D\) \emph{arc-faithful}
if the over- and under-arc meridians at each crossing are distinct in
\(\pi_1(M)\). For complete holonomy this is equivalent to \(D\)-smoothness,
meaning that a non-pinched octahedral coloring realizes the holonomy
up to conjugacy
\cite[Definition~26, Theorem~3]{McPhailSnyder}.

For reduced alternating diagrams of hyperbolic knots, the corresponding
non-peripherality of the four-term decomposition was established by
Garoufalidis--Moffatt--Thurston and independently by Sakuma--Yokota
\cite{GMT,SakumaYokota}. McPhail--Snyder also records existence for
hyperbolic knots with at most twelve crossings, based on the
boundary-parabolic representation data of Kim and Choi \cite{KimChoi},
and asks whether
every hyperbolic knot has an arc-faithful diagram
\cite[\S5, Conjecture~29]{McPhailSnyder}.

\begin{theorem}\label{thm:main}
Let \(K\subset S^3\) be a hyperbolic knot and let \(\beta\in B_n\) be any
braid whose closure is \(K\). Each word for \(\beta\) can be replaced by a
word for the same element of \(B_n\) whose closed braid diagram is
arc-faithful. The complete holonomy is realized by a non-pinched
octahedral coloring on this diagram.
\end{theorem}

Alexander's theorem supplies a closed braid for every oriented knot
\cite{Alexander}, so Theorem~\ref{thm:main} proves the stated conjecture.

\begin{samepage}
The strand meridians at each transverse slice generate the knot group.
If their parabolic images all shared a fixed point, moving that point to
infinity would make every image a translation. The quotient of hyperbolic
space by such a group has infinite volume. Thus at a pinched crossing,
where the incoming pair shares a fixed point, the finite volume of \(M\)
guarantees a third strand with a different fixed point.
\par
\end{samepage}

Once that strand is adjacent to the pair, the three-strand identity
\(a=s^{-1}asa s^{-1}\) for adjacent Artin generators \(a,s\)
replaces the pinched crossing by five non-pinched crossings. A strand
farther away can be brought to the pair using crossings supported outside
the pair. The construction preserves the braid axis, strand count, and
writhe.

\section{Parabolic colors on a braid}\label{sec:colors}

Let \(B_n\) have Artin generators \(\sigma_1,\dots,\sigma_{n-1}\)
\cite{Artin}. Read words in time order from left to right. In our convention
\(\sigma_i\) has the right incoming strand passing over the left.
The Wirtinger relations give the Hurwitz rules on the adjacent colors
\((x,y)\) in positions \(i,i+1\).
\begin{equation}\label{eq:hurwitz}
(x,y)\cdot \sigma_i=(y,y^{-1}xy),\qquad
(x,y)\cdot \sigma_i^{-1}=(xyx^{-1},x).
\end{equation}
These mutually inverse rules define a right action of \(B_n\) on tuples
of elements of any group. Indeed the rules satisfy the Artin braid
relations, and operations on disjoint pairs commute.

Given a representation \(\rho\) of the closed braid link group, we color
each strand segment by the image of its Wirtinger meridian. At a regular
slice the colors are an ordered \(n\)-tuple. The meridians at the initial
slice generate the link group. Each new meridian is obtained by
propagating the tuple using \eqref{eq:hurwitz}, and closing the braid
imposes the relations identifying the final tuple with the initial one.
The subgroup generated by the colors is constant from slice to slice,
since either pair in \eqref{eq:hurwitz} generates the same subgroup as its
input pair. It follows that every slice generates the image of \(\rho\).

For a nontrivial parabolic \(x\in\PSL\), write \(\fix(x)\) for its unique
fixed point in \(\CP\). Conjugation transports this point according to
\(\fix(gxg^{-1})=g\fix(x)\), where \(g\) acts by a M\"obius transformation.
Writing \(P=\fix(x)\) and \(Q=\fix(y)\), the ordered fixed points of the
outputs in \eqref{eq:hurwitz} are
\begin{equation}\label{eq:point-update}
(Q,y^{-1}P)\quad\text{and}\quad (xQ,P),
\end{equation}
respectively. Since \(yQ=Q\) and \(xP=P\), distinct incoming points give
distinct outgoing points. Consequently a sequence of non-pinched
crossings stays non-pinched when reversed from its output tuple.
If all colors at a slice fix the same point, conjugating that point
to infinity makes them translations \(z\mapsto z+c\), so they generate
an abelian group.

\section{Eliminating pinched crossings}\label{sec:replacement}

\begin{lemma}[Adjacent third strand]\label{lem:local}
Suppose the colors \(x,y\) at a crossing \(a^\varepsilon\), with
\(a=\sigma_i\) and \(\varepsilon\in\{1,-1\}\), fix the same point \(P\).
Suppose a neighboring third strand has color \(h\) with
\(\fix(h)=Q\ne P\). Let \(s\) be the generator between the third strand
and the nearer strand of the pair. Then the word
\begin{equation}\label{eq:replacement}
W_\varepsilon=s^{-\varepsilon}a^\varepsilon s^\varepsilon
                  a^\varepsilon s^{-\varepsilon}
\end{equation}
represents \(a^\varepsilon\), gives the same output color tuple,
and has five non-pinched crossings.
\end{lemma}

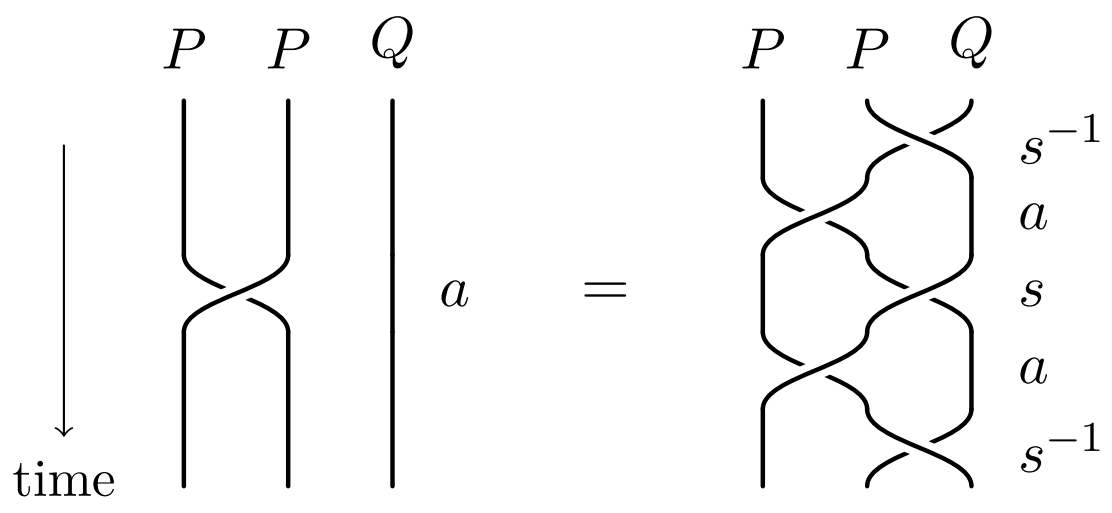
\begin{figure}[tbp]
\centering
\begin{tikzpicture}[x=.69cm,y=.51cm,line width=.84pt,line cap=round]
  \foreach \pos/\point in {0/P,1/P,2/Q}{
    \node[above=2pt] at (\pos,0) {\(\point\)};
    \draw (\pos,0) -- (\pos,-2);
    \draw (\pos,-3) -- (\pos,-5);
  }
  \braidlayer{0}{-2}{1}
  \node[anchor=west] at (2.24,-2.5) {\(a\)};
  \draw[->,thin] (-1.15,-.58) -- (-1.15,-4.36);
  \node[below,font=\small] at (-1.15,-4.41) {time};
  \node[font=\large] at (4.04,-2.5) {\(=\)};
  \begin{scope}[shift={(5.55,0)}]
    \foreach \pos/\point in {0/P,1/P,2/Q}
      \node[above=2pt] at (\pos,0) {\(\point\)};
    \braidlayer{1}{0}{-1}
    \braidlayer{0}{-1}{1}
    \braidlayer{1}{-2}{1}
    \braidlayer{0}{-3}{1}
    \braidlayer{1}{-4}{-1}
    \foreach \level/\letter in {-.5/{s^{-1}},-1.5/a,-2.5/s,-3.5/a,-4.5/{s^{-1}}}
      \node[anchor=west] at (2.24,\level) {\(\letter\)};
  \end{scope}
\end{tikzpicture}
\caption{The positive replacement with the third strand on the right.
The labels \(P,P,Q\) are the incoming fixed points. The single crossing
on the left is pinched. All five crossings on the right are non-pinched
when \(P\ne Q\).}
\label{fig:replacement}
\end{figure}

\begin{proof}
The Artin relation \(asa=sas\) gives \(W_1=a\), and inversion gives
\(W_{-1}=a^{-1}\). Figure~\ref{fig:replacement} shows the positive case
with the third strand on the right. Equality in the braid group also
gives equality of the Hurwitz actions, including the output on the third
strand.

\begin{samepage}
For the positive case, the input tuple is \((x,y,h)\) when the third
strand is on the right and \((h,x,y)\) when it is on the left.
Set \(u=yhy^{-1}\). The following table records the two incoming
points at each letter of \(W_1\), in their strand order.
\[
\begin{array}{@{}c@{\qquad}c@{\qquad}c@{}}
\toprule
\text{Letter} & \text{Third strand on the right} & \text{Third strand on the left}\\
\midrule
s^{-1} & (P,Q) & (Q,P)\\
a      & (P,yQ) & (Q,P)\\
s      & (u^{-1}P,P) & (hP,P)\\
a      & (yQ,P) & (y^{-1}hP,y^{-1}Q)\\
s^{-1} & (Q,h^{-1}P) & (P,y^{-1}Q)\\
\bottomrule
\end{array}
\]
\par
\end{samepage}
On the right, the first letter moves \(u\) into the middle, and the
second moves it to the left while conjugating \(x\) by \(u^{-1}\).
The third then gives the rightmost point
\(y^{-1}u^{-1}P=h^{-1}P\). The fourth moves \(h\) back to the middle,
as \(y^{-1}uy=h\). On the left, the first two letters give the tuple
\((hxh^{-1},y,y^{-1}hy)\). The third conjugates \(hxh^{-1}\) by
\(y^{-1}\) and the fourth gives rightmost color \(y^{-1}xy\).
These updates give the entries of the table by \eqref{eq:point-update}.

Since \(yP=P\), both \(yQ\) and \(y^{-1}Q\) differ from \(P\).
The point \(u^{-1}P\) differs from \(P\) because \(u\) fixes only \(yQ\).
Also \(hP\ne P\) by uniqueness of \(Q\), and \(h^{-1}P\ne Q\)
since \(hQ=Q\). Finally \(hP\ne hQ=Q\) implies
\(y^{-1}hP\ne y^{-1}Q\). Thus the points differ in every row.

For the negative case, start the positive construction at the output of
\(a^{-1}\). The two colors of that pair still have fixed point \(P\)
and the third color is unchanged. The positive word returns to the
original input of \(a^{-1}\). Reverse this path, using the reversibility
of non-pinched crossings, to obtain \(W_{-1}\).
\end{proof}

\begin{theorem}[Rewriting with a prescribed coloring]\label{thm:rewrite}
Let \(w\) be a braid word on \(n\ge 2\) strands and let its input tuple
consist of nontrivial parabolics in \(\PSL\) with no common fixed point.
There is a word \(w'\) for the same braid element such that every crossing
of \(w'\), colored from that tuple, is non-pinched. If \(w\) has length
\(L\) and \(p\) pinched crossings, then
\begin{equation}\label{eq:length}
|w'|\le L+(2n-2)p.
\end{equation}
\end{theorem}

\begin{proof}
At a pinched crossing \(a^\varepsilon\) the two participating colors
fix a point \(P\). The tuple at that slice generates the same subgroup as
the input tuple, so some third color has a fixed point different from
\(P\). Choose such a strand closest to the pair, with \(d\) intervening
strands. All intervening colors fix \(P\).

\begin{samepage}
Bring the chosen strand to the pair using \(d\) crossings. To move its
current color \(h\) left past a color \(g\) that fixes \(P\), use
\((g,h)\mapsto(ghg^{-1},g)\). To move it right, use
\((h,g)\mapsto(g,g^{-1}hg)\). In either case \(g\) is unchanged and
the moving color retains a fixed point different from \(P\), since
both \(g\) and \(g^{-1}\) fix \(P\). Thus the transport word \(C\) is
non-pinched. Its letters are supported outside the original pair,
so \(C\) commutes with \(a\).
\par
\end{samepage}

Replace \(a^\varepsilon\) by \(CW_\varepsilon C^{-1}\), using
Lemma~\ref{lem:local} at the new position of the third strand. The middle
word acts as \(a^\varepsilon\) and leaves the colors outside the original
pair unchanged. The return \(C^{-1}\) therefore reverses the transport
on exactly the same outside colors and remains non-pinched. The new
block represents \(Ca^\varepsilon C^{-1}=a^\varepsilon\) and gives
the original output tuple.

At every replaced letter this output agreement preserves the input
tuple of the next original letter. Replacing the \(p\) pinched crossings
in order therefore gives a word with every crossing non-pinched.
For \(n\ge 3\) each block has length \(5+2d\le 2n-1\), since
\(d\le n-3\). The length increases by at most \(2n-2\) at each
replacement. For \(n=2\) the hypothesis and the invariance of the
generated subgroup imply \(p=0\).
\end{proof}

Theorem~\ref{thm:rewrite} applies to the meridian coloring of any
oriented link for a prescribed representation with the stated
parabolic and fixed-point hypotheses.

\section{Hyperbolic holonomy and the closed diagram}\label{sec:holonomy}

\begin{proof}[Proof of Theorem~\ref{thm:main}]
The complete holonomy images of the slice meridians are nontrivial
parabolics generating a finite-covolume Kleinian group. A common fixed
point would make this group a discrete translation group. In the
upper half-space model its fundamental domain would have constant
positive Euclidean cross-sectional area at every height \(t>0\), and
thus infinite volume by the density \(t^{-3}\).

Apply Theorem~\ref{thm:rewrite}. The isotopy of the two braids relative
to their endpoints identifies the initial meridians and transports
the complete holonomy. The distinct fixed points at every crossing
imply distinct over- and under-meridians in the knot group, giving
arc-faithfulness. The complete holonomy lifts to \(\SL\), with its
decoration given by the unique invariant lines of its parabolic
meridians. McPhail--Snyder's fixed-point criterion gives the
non-pinched octahedral coloring \cite[Corollary~25]{McPhailSnyder}.
\end{proof}

The exponent sum factors through \(B_n\), so the construction
preserves the writhe and the associated diagram framing. Starting at
minimal braid index gives an arc-faithful diagram at that index.

\begin{samepage}
\subsection{A stabilized figure-eight braid}
Let \(b=(\sigma_1\sigma_2^{-1})^2\) be the standard alternating
figure-eight braid \cite[\S4.2]{MOV}. Consider the positive stabilization
\(b\sigma_3\in B_4\). If \((x_1,x_2,x_3)\) is the complete holonomy tuple
for \(b\), the stabilization repeats its third meridian, giving the
input \((x_1,x_2,x_3,x_3)\). The final \(\sigma_3\) is pinched.
The preceding \(\sigma_2^{-1}\) is non-pinched by arc-faithfulness
of the alternating diagram, so its outputs \(x_2,x_3\) have distinct
fixed points. Using the strand in position \(2\) in
Lemma~\ref{lem:local} gives a nine-crossing word for \(b\sigma_3\)
with all crossings non-pinched.
\[
(\sigma_1\sigma_2^{-1})^2\,
\sigma_2^{-1}\sigma_3\sigma_2\sigma_3\sigma_2^{-1}.
\]
\par
\end{samepage}

\subsection{Distinct regions at each crossing}

\begin{proposition}\label{prop:regions}
The closed braid in Theorem~\ref{thm:main} can be chosen so that at each
crossing the four incident sectors belong to four distinct complementary
regions. This requires at most \(2(n-1)\) additional letters.
\end{proposition}

\begin{proof}
Use the standard annular projection of a braid, with crossings at distinct
angles around the axis. Number the radial channels \(0,\ldots,n\) by
the number of strand positions closer to the inner boundary. This
number is locally constant off the diagram. At a crossing of index
\(i\), one sector lies in each of channels \(i-1\) and \(i+1\).
The other two sectors lie in channel \(i\) just before and just after
the crossing. Channel \(i\) is cut into components by the angular
positions of the crossings with index \(i\). If that index occurs at
least twice, the two sectors in channel \(i\) are in distinct
components. This also describes the complementary regions in the
planar closure, with channels \(0\) and \(n\) extending to the two
regions outside the annulus.

Each generator index occurs at least once, since otherwise the braid
permutation preserves the two nonempty blocks separated by that index,
contradicting that the closure is a knot. If an index occurs exactly
once, replace its letter \(\sigma_i^\varepsilon\) by
\(\sigma_i^\varepsilon\sigma_i^{-\varepsilon}\sigma_i^\varepsilon\).
All three crossings are non-pinched by reversibility, and their product
is the original letter. This adds two letters for each affected index.
\end{proof}

\section{Consequences for diagram potentials}\label{sec:consequences}

By McPhail--Snyder's correspondence, the non-pinched coloring realizes
the complete holonomy by a solution of the four-term gluing equations
with shapes in \(\mathbb C\setminus\{0,1\}\). The edges of this
decomposition are homotopically non-peripheral
\cite[Corollary~27(a)]{McPhailSnyder}. Here the decomposition retains the
two auxiliary ideal vertices. Nondegeneracy of a shape means that the
four ideal vertices of its tetrahedron are distinct. This condition
allows real shapes and either orientation of straightened tetrahedra.

For a hyperbolic knot the Volume Conjecture predicts
\begin{equation}\label{eq:volume}
\lim_{N\to\infty}\frac{2\pi}{N}\log\bigl|\langle K\rangle_N\bigr|
=\Vol(M),
\end{equation}
where \(\langle K\rangle_N\) is the Kashaev invariant, equivalently the
reduced \(N\)-colored Jones polynomial at \(\exp(2\pi i/N)\) in the
usual zero-framing convention \cite{Kashaev,MurakamiMurakami}.
The diagram potentials used in saddle-point approaches are sums of
dilogarithms evaluated at shape parameters \cite{Yokota,McPhailSnyder}.
At a pinched crossing a shape equal to \(1\) lies at a branch point
of the dilogarithm. Theorem~\ref{thm:main} supplies, for every hyperbolic
knot, a diagram in which the complete holonomy avoids this obstruction.

\begin{samepage}
Let \(\mu\) be a logarithmic meridian parameter of the decorated lift.
In logarithmic segment coordinates \(\beta_j\), the non-pinched
complete point \(\boldsymbol\beta_*\) satisfies the generalized
critical equations of \(\Phi_{D,\mu}\) \cite[\S5]{McPhailSnyder}.
\[
\exp\left(\frac{\partial\Phi_{D,\mu}}{\partial\beta_j}
          (\boldsymbol\beta_*)\right)=1
\qquad\text{for every segment }j.
\]
Local branches of the logarithms and dilogarithms give a holomorphic
germ at \(\boldsymbol\beta_*\). Subtracting the linear term of this
germ makes it stationary there and preserves the exponentiated
equations, since its first derivatives at \(\boldsymbol\beta_*\)
lie in \(2\pi i\mathbb Z\).
\par
\end{samepage}

For the state integrals of the geometrically twisted invariant
\cite[\S2.2]{McPhailSnyderState}, one may choose a gauge in which the
shapes avoid the unit circle \cite[Corollary~21]{McPhailSnyder}.
Opening a segment incident to the unbounded region produces the
cut diagram with the same crossing data.

The analytic conclusion is the existence of a smooth local phase at
the complete holonomy. Deriving \eqref{eq:volume} from a quantum state
sum also requires identifying its additive phase and integration cycle,
showing that the complete point contributes, and controlling other
contributions and cancellations.

\paragraph{Tool Use.}
OpenAI Codex helped develop and check the proof, locate references,
and prepare the \LaTeX{} manuscript and TikZ braid diagram.


\begin{thebibliography}{99}
\small
\raggedright
\setlength{\itemsep}{3pt}
\bibitem{Alexander}
J.~W. Alexander,
\emph{A lemma on systems of knotted curves},
Proc. Natl. Acad. Sci. USA \textbf{9} (1923), 93--95.
\href{https://doi.org/10.1073/pnas.9.3.93}{DOI}.

\bibitem{Artin}
E.~Artin,
\emph{Theory of braids},
Ann. of Math. (2) \textbf{48} (1947), 101--126.
\href{https://doi.org/10.2307/1969218}{DOI}.

\bibitem{GMT}
S.~Garoufalidis, I.~Moffatt, and D.~P. Thurston,
\emph{Non-peripheral ideal decompositions of alternating knots},
\href{https://arxiv.org/abs/1610.09901}{arXiv:1610.09901}, 2016.

\bibitem{Kashaev}
R.~M. Kashaev,
\emph{The hyperbolic volume of knots from the quantum dilogarithm},
Lett. Math. Phys. \textbf{39} (1997), 269--275.
\href{https://doi.org/10.1023/A:1007364912784}{DOI}.

\bibitem{KKY}
H.~Kim, S.~Kim, and S.~Yoon,
\emph{Octahedral developing of knot complement. I:
Pseudo-hyperbolic structure},
Geom. Dedicata \textbf{197} (2018), 123--172.
\href{https://doi.org/10.1007/s10711-018-0323-8}{DOI}.

\bibitem{KimChoi}
S.~Kim and P.~Choi,
\emph{Diagram.site: Boundary parabolic representations for knots},
\href{https://diagram.site/}{diagram.site}.

\bibitem{McPhailSnyder}
C.~McPhail-Snyder,
\emph{Octahedral coordinates from the Wirtinger presentation},
Geom. Dedicata \textbf{219} (2025), article~53.
\href{https://doi.org/10.1007/s10711-025-01012-7}{DOI}.
Statement numbers in the text follow
\href{https://arxiv.org/html/2404.19155v2}{arXiv:2404.19155v2}.

\bibitem{McPhailSnyderState}
C.~McPhail-Snyder,
\emph{State integrals for the quantized \(\SL\) Chern--Simons invariant},
\href{https://arxiv.org/abs/2601.05136v2}{arXiv:2601.05136v2}, 2026.

\bibitem{MOV}
S.~Morier-Genoud, V.~Ovsienko, and A.~P. Veselov,
\emph{Burau representation of braid groups and \(q\)-rationals},
Int. Math. Res. Not. IMRN (2024), no.~10, 8618--8627.
\href{https://doi.org/10.1093/imrn/rnad318}{DOI}.

\bibitem{MurakamiMurakami}
H.~Murakami and J.~Murakami,
\emph{The colored Jones polynomials and the simplicial volume of a knot},
Acta Math. \textbf{186} (2001), 85--104.
\href{https://doi.org/10.1007/BF02392716}{DOI}.

\bibitem{SakumaYokota}
M.~Sakuma and Y.~Yokota,
\emph{An application of non-positively curved cubings of alternating links},
Proc. Amer. Math. Soc. \textbf{146} (2018), 3167--3178.
\href{https://doi.org/10.1090/proc/13918}{DOI}.

\bibitem{Yokota}
Y.~Yokota,
\emph{On the potential functions for the hyperbolic structures of a knot
complement},
Geom. Topol. Monogr. \textbf{4} (2002), 303--311.
\href{https://doi.org/10.2140/gtm.2002.4.303}{DOI}.
\end{thebibliography}
\end{document}